\documentclass[10pt, a4paper, reqno]{amsart}
\usepackage{amscd}
\usepackage{amsmath}
\usepackage{amstext}
\usepackage{amsthm}
\usepackage{amssymb}
\usepackage{amsopn}
\usepackage{amssymb}
\usepackage{amsxtra}
\usepackage{amsfonts}
\usepackage{fancyhdr}
\usepackage{paralist}
\usepackage[mathcal]{euscript}
\usepackage[english]{babel}
\usepackage[latin1]{inputenc}

\usepackage{tikz, bbding, subfigure, pifont, color, manfnt}

\newcommand{\gray}[1]{\textcolor{gray}{#1}}

\newcommand{\KreisKleinSchwarz}{\ensuremath{{\scriptscriptstyle\text{\ding{108}}}}}
\newcommand{\KreisKleinWeiss}{\ensuremath{\text{\begin{picture}(0,0)(0,0)\put(2,1.64){\circle{2.5}}\end{picture}}}}
\newcommand{\KreisKleinGrau}{\ensuremath{\gray{\scriptscriptstyle\text{\ding{108}}}}}

\newcommand{\KreisGrossSchwarz}{\ensuremath{{\scriptstyle\text{\ding{108}}}}}
\newcommand{\KreisGrossWeiss}{\ensuremath{\text{\begin{picture}(8,0)(0,0)\put(3.5,2.25){\circle{4.5}}\end{picture}}}}
\newcommand{\KreisGrossGrau}{\ensuremath{\gray{\scriptstyle\text{\ding{108}}}}}

\newcommand{\QuadKleinSchwarz}{\ensuremath{{\scriptscriptstyle\blacksquare}}}
\newcommand{\QuadKleinWeiss}{\ensuremath{{\scriptscriptstyle\square}}}
\newcommand{\QuadKleinGrau}{\ensuremath{\gray{\scriptscriptstyle\blacksquare}}}

\newcommand{\QuadGrossSchwarz}{\ensuremath{{\scriptstyle\blacksquare}}}
\newcommand{\QuadGrossWeiss}{\ensuremath{{\scriptstyle\square}}}
\newcommand{\QuadGrossGrau}{\ensuremath{\gray{\scriptstyle\blacksquare}}}

\newcommand{\DreieckKleinSchwarz}{\ensuremath{{\scriptscriptstyle\text{\ding{115}}}}}
\newcommand{\DreieckKleinWeiss}{\ensuremath{{\scriptscriptstyle\bigtriangleup}}}
\newcommand{\DreieckKleinGrau}{\ensuremath{\gray{\scriptscriptstyle\text{\ding{115}}}}}

\newcommand{\DreieckGrossSchwarz}{\ensuremath{{\scriptstyle\text{\ding{115}}}}}
\newcommand{\DreieckGrossWeiss}{\ensuremath{{\scriptstyle\bigtriangleup}}}
\newcommand{\DreieckGrossGrau}{\ensuremath{\gray{\scriptstyle\text{\ding{115}}}}}

\newtheoremstyle{normal}% normale Schrift
{5pt}% hSpace abovei
{5pt}% hSpace belowi
{\normalfont}% hBody fonti
{}% hIndent amounti1
{\bfseries}% hTheorem head fonti
{}% Punctuation after theorem headi
{0.4em}% hSpace after theorem headi2
{\bfseries{\thmname{#1}\thmnumber{ #2}.\thmnote{ \hspace{0.5em}(#3)\newline}}}% hTheorem head spec (can be left empty, meaning `normal')

\newtheoremstyle{kursiv}% kursive Schrift
{5pt}% hSpace abovei
{5pt}% hSpace belowi
{\itshape}% hBody fonti
{}% hIndent amounti1
{\bfseries}% hTheorem head fonti
{}% Punctuation after theorem headi
{0.4em}% hSpace after theorem headi2
{\bfseries{\thmname{#1}\thmnumber{ #2}.\thmnote{ \hspace{0.5em}(#3)\newline}}}% hTheorem head spec (can be left empty, meaning `normal')

\theoremstyle{kursiv}
\newtheorem{thm}{Theorem}

\newtheorem{dfn}[thm]{Definition}

\newtheorem{exer}[thm]{Exercise}
\newtheorem{ques}[thm]{Question}
\newtheorem{rem}[thm]{Remark}
\newtheorem{sol}[thm]{Solution}

\renewcommand{\emptyset}{\varnothing}
\renewcommand{\epsilon}{\varepsilon}

\begin{document}

\title{A workshop for high school students\\on naive set theory}

\author{Sven-Ake Wegner\,\MakeLowercase{$^{\text{a}}$}\vspace{25pt}}

\renewcommand{\thefootnote}{}

\hspace{-1000pt}\footnote{\hspace{5.5pt}2010 \emph{Mathematics Subject Classification}: Primary 97D80; Secondary 97C90, 97D40.}
\hspace{-1000pt}\footnote{\hspace{5.5pt}\emph{Key words and phrases}: workshop for high school students, naive set theory, student life cyle.\vspace{1.6pt}}
\hspace{-1000pt}\footnote{$^{\text{a}}$\,Bergische Universit\"at Wuppertal, FB C -- Mathematik und Naturwissenschaften, Gau\ss{}stra\ss{}e 20, D-42119 Wuppertal,\linebreak\phantom{x}\hspace{12.5pt}Germany, Phone:\hspace{1.2pt}\hspace{1.2pt}+49\hspace{1.2pt}(0)\hspace{1.2pt}202\hspace{1.2pt}/\hspace{1.2pt}439\hspace{1.2pt}-\hspace{1.2pt}2531, Fax:\hspace{1.2pt}\hspace{1.2pt}+49\hspace{1.2pt}(0)\hspace{1.2pt}202\hspace{1.2pt}/\hspace{1.2pt}439\hspace{1.2pt}-\hspace{1.2pt}3724, E-Mail: wegner@math.uni-wuppertal.de.}

\begin{abstract}
In this article we present the prototype of a workshop on naive set theory designed for high school students in or around the seventh year of primary education. Our concept is based on two events which the author organized in 2006 and 2010 for students of elementary school and high school, respectively. The  article also includes a practice report on the two workshops.
\end{abstract}

\maketitle

\begin{picture}(0,0)
\put(195,95){\small April 16, 2014}
\end{picture}

\vspace{-20pt}

\section{Introduction}\label{SEC-1}

Events for high school students or students of elementary school in altering configurations and with different aims are a regular component of modern university teaching portfolio. The importance of events that supplement regular classes on all levels of school education increases constantly. In order to support and encourage interested students and to attract prospective university students the universities and the schools organize events such as open days, information sessions, workshops etc.
\smallskip
\\When planning an event of the above type, it is a natural question, which topic should be discussed and which teaching methods should come into play. In the case of a workshop for high school students at the university the topic should give an insight into academic life. In addition, the university teaching and learning concepts should be applied during the workshop in order to give the participants also an insight in their prospective student life. In the case of a course at high or elementary school, the topics covered by the usual curriculum should be supplemented and the participants should be challenged according to their personal level of performance---possibly in contrast to regular classes. In order to meet the above objectives, a workshop must include exercise sessions in which the students work on questions individually or in small groups. The level of support during these exercise sessions must vary from intensive discussions to no support at all. Exercises form a substantial part of mathematical studies, in particular homework which has to be done without any assistance from a teacher is indispensable for mathematics education. The choice of a suitable topic, on the one hand challenging and on the other hand not of unreachable difficulty, in combination with the design of a workshop concept with the above features is a non-trivial task for the designated organizer.
\smallskip
\\This paper is based on two different events which the author organized in 2006 and 2010: In Section 4 we report on a workshop for students of elementary school in their third or fourth year of primary education which was part of a program initiated by a regional group of the German Association for Highly Talented Children (DGhK-OWL). The second event, on which we report in Section 4, was a workshop for high school students in their seventh or eighth year of primary education. It was carried out in the context of the Girl's Day 2010 and took place at the University of Paderborn. In Section 3 we present the prototype of a workshop on naive set theory based on the aforementioned events. Before this, we outline our motivation for this topic in the context of student workshops in Section 2 below.
\smallskip
\\We stress that all statements about execution and success of the workshops we organized are based on the personal assessment of the author; a formal evaluation was not carried out. We understand the current paper as a \textquotedblleft{}practice report\textquotedblright{}. Our aim is to share our experiences with prospective organizers of future workshops.

\section{Why Set Theory?}

In view of the aims mentioned in Section 1, it is not easy to find a topic which allows for the design of a workshop that meets all our objectives. Looking for a topic with the adequate level of difficulty and an appropriate distance to regular mathematics lessons, the author organized his first workshop on naive set theory in 2006, see Section 4. Below we list arguments for the latter topic as a possible subject of student workshops.

\begin{compactitem}

\vspace{5pt}

\item[1.] Naive set theory is often hidden but ubiquitous in mathematics. Already on school level, naive set theory is used as a language (set of natural numbers, solution set, zero set, primitive integral etc.) and this language is indispensable for formulating complex mathematical statements. In this sense naive set theory is relevant for all branches of mathematics and all branches of science that employ mathematical techniques.

\vspace{3pt}

\item[2.] The study of set theory can be used to strengthen the ability of formal and abstract reasoning. When dealing with mathematically interested and talented students, it seems reasonable to promote these skills in order to complement the regular courses which sometimes focus more on applications of mathematics in business or technology.

\vspace{3pt}

\item[3.] Gaining knowledge on set theory is a preparation for mathematical studies or studies in any related field. Indeed, a short introduction to naive set theory can be found at the beginning of many first year university courses or textbooks.

\vspace{3pt}

\item[4.] The study of naive set theory requires no special previous mathematics education; only an appropriate repertory of examples is needed. Therefore, a workshop on naive set theory can be adjusted to a large scale of possible participants according to age, education level etc.

\vspace{3pt}

\item[5.] In the context of naive set theory it is possible to reach interesting problems or applications in relatively short time, e.g., Russell's antinomy, see for instance Potter \cite[Section 2.3]{Potter}, or the formal definition of a map. The latter illustrates the benefit of the notion of a set in the context of topics that arise in regular lessons.  Russell's antinomy, on the other hand, illustrates the limitations of naive set theory.

\vspace{3pt}

\item[6.] Naive set theory allows for the use of different teaching methods, many exercises can be completed in small groups or individually, see Narli, Baser \cite{NarliBaser}. In particular the first exercises, operated with \textquotedblleft{}touchable\textquotedblright{} sets are suitable to be solved via teamwork in a natural way.

\end{compactitem}\vspace{5pt}

We point out that the above statements refer to the specific situation of optional and supplementary events designed for high school students which are interested in mathematics. We distance ourself from a treatment of set theory within regular mathematics school education. In particular, we do not support the 1960s philosophy of so-called \textquotedblleft{}New Math\textquotedblright{}, cf.~Kline \cite{Kline}, and reject the philosophy that teaching mathematics at elementary school should start with set theory. We stress that the workshop proposed in Section 3 and the school education proposed in textbooks featuring the \textquotedblleft{}New Math\textquotedblright{} are substantially different in their intention, design and context.

\section{The Workshop}

In the sequel, we present the prototype of a workshop on naive set theory. For variants concerning teaching methods we refer to the relevant didactics literature. More information on naive set theory can be found at the beginning of many undergraduate textbooks, e.g., Gaughan \cite{Gaughan} or Rosenlicht \cite{Rosenlicht}. For more details we refer to Halmos \cite{Halmos}. Above, we distanced ourself from the philosophy of \textquotedblleft{}New Math\textquotedblright{}. However, we note that the teaching material developed in the 1960s, e.g., \cite{Matema}, can be useful for the preparation of a workshop on naive set theory. The design of a concrete workshop depends on the underlying conditions (available time, number of participants, age and educational background, etc.). In the sequel we assume the following.
\smallskip
\\\textbf{Type of Event:} Supplementary workshop for mathematically interested or talented high school students. 
\smallskip
\\\textbf{Age:} Around the seventh year of primary education. The variation within the group should not be more than one year. Under special conditions it is also possible to run the workshop with students of elementary school, cf.~Section 4.
\smallskip
\\\textbf{Number of Participants:} Less or equal 15.
\smallskip
\\\textbf{Time:} 2--3 hours.
\smallskip
\\\textbf{Classroom:} Classroom with blackboard and the possibility for team work sections.
\smallskip
\\\textbf{Lecturers:} One main lecturer and, depending on the number of participants, additional lecturers for the team work sections. 
\smallskip
\\\textbf{Materials:} For every team a set of 48 shapes made for instance of plastic, wood or cardboard with the following properties in all possible combinations.

\begin{center}
\begin{tabular}{rl}
colors:& red, yellow, green, blue,\\
shapes:& circular, triangular, quadratic,\\
sizes:& big, small,\\
surfaces:& plain, abrasive,\\
\end{tabular}
\end{center}

one big sheet of paper (approximately 40$\times$60cm), one board marker, notepad and pencil for every participant.
\smallskip
\\In the two workshops discussed in Section 4, the author used the so-called \textquotedblleft{}Matema Begriffsspiel\textquotedblright{} \cite{Matema} for the intuitive approach below. However, it is easy to adjust the outline below according to other sets of shapes. In the sequel, we note the main steps of our workshop. We enrich our description with concrete tasks and solutions. Depending on the age of the participants, their previous education and skills, the sections below can be extended or shortened.

\subsection*{I.~Intuitive approach} To begin with, the participants develop an intuitive idea of sets by working on the following exercise.

\begin{exer}
Sort the shapes.
\end{exer}

During the work on the exercise, the participants get acquainted with each other and with the teaching material (shapes). In small teams they discuss various options for sorting the shapes. Later they discuss their conclusions together with the lecturer. The results should include the following.

\begin{compactitem}\vspace{5pt}

\item[(a)] The task in Exercise 1 is not well-posed. Each of the shapes has several properties and there is no natural principle for a sorting as it is for instance the case in other typical \textquotedblleft{}sorting tasks\textquotedblright{} (coins can be sorted according to their value, words can be sorted lexicographically, etc.).

\vspace{3pt}

\item[(b)] It is thus necessary to decide for a sorting principle. For instance, one can sort according to the four colors and ignore the other properties to get four piles of shapes. But one can also sort with respect to sizes. Then one gets two piles of shapes, etc.

\vspace{3pt}

\item[(c)] If one proceeds as in (b), then the shapes on every pile (for instance all yellow shapes) remain \emph{distinguishable} via there other properties (shape, surface and size). At the same time they are \emph{distinct} (in this case by their yellow color).

\vspace{5pt}

\end{compactitem}

Concerning (a), the lecturer can assist during the discussion by giving the other examples for sorting problems (coins and words) mentioned above. Concerning (b), a concretization of the task in the exercise can help, for instance the lecturer can ask for a sorting according to colors. In addition, the lecturer can then ask to compare this concrete task with the (non-concrete) original exercise. It cannot be expected that the participants will find the conclusions mentioned in (c) in this form. However, the lecturer can check whether they developed an intuitive idea of sets. If, for instance, the participants sort according to colors, the lecturer can ask why a particular shape cannot belong to several piles at the same time or to no pile at all.

\subsection*{II. Naive Definition} In a moderated discussion the main points of the definition of a set are summarized, see (c) above. The lecturer notes the definition on the blackboard, for instance in the following form, cf.~Dauben \cite[p.~170]{Dauben}.

\medskip

\begin{dfn} By a \emph{set} we mean any collection into a whole of definite, distinct objects of our perception or of our thought.
\end{dfn}

This definition is \emph{naive}; the words \textquotedblleft{}collection\textquotedblright{}, \textquotedblleft{}definite\textquotedblright{},  \textquotedblleft{}distinct\textquotedblright{}, \textquotedblleft{}object\textquotedblright{} are not defined in a formal way. The latter is not possible in our context and in particular \emph{not intended}. The above formulation has to be illustrated by considering concrete instances of sets. It should be regarded from this more practical rather than from a philosophical point of view. After giving the above definition, the lecturer explains the following notation.

\begin{dfn} If an object belongs to a set, then we say that it is an \emph{element} of the set. In order to describe a set, we may \emph{list} its elements or we may \emph{characterize} them by specifying a certain rule which regulates which objects are in the set and which are not, e.g., $\{{\scriptscriptstyle\square},\,{\scriptstyle\square},\,{\scriptscriptstyle\bigtriangleup}\}$ or $\{P\:;\:P\text{ is blue}\}$. In this case we write ${\scriptstyle\square}\in\{{\scriptscriptstyle\square},\,{\scriptstyle\square},\,{\scriptscriptstyle\bigtriangleup}\}$. On the other hand we write ${\begin{picture}(7,0)(0,0)\put(3.5,2.5){\circle{4.5}}\end{picture}}\not\in\{{\scriptscriptstyle\square},\,{\scriptstyle\square},\,{\scriptscriptstyle\bigtriangleup}\}$, if an object does not belong to a set. To start with, we restrict to the set of all shapes which we denote by $U$. The sets above are \emph{subsets} of $U$, we write $\{P\:;\:P\text{ is blue}\}\subseteq U$.
\end{dfn}

\subsection*{III.~Operations with Sets} After fixing the notation, we introduce the set operations (intersection, union, etc.). Similar to the notion of a set itself, the latter can first be discussed on an intuitive level using the shapes via the following exercise. The lecturer has to take care about the \emph{process} in which the sets below are formed and leads the discussion in the direction of set operations.

\begin{exer} Form the sets of all shapes, which are
\smallskip
\\\begin{tabular}{llllll}\vspace{1pt}
(1) & big and blue, & & (4) & circular and quadratic, \\\vspace{1pt}
(2) & plain or abrasive, & & (5) & red or blue or yellow,\\\vspace{1pt}
(3) & not abrasive, & & (6) & circular and plain and not green.
\end{tabular}
\end{exer}

Depending on the participant's competencies, not all of the above examples have to be treated via shapes. Beginning with the first example, Venn diagrams are drawn under the guidance of the lecturer, cf.~Figure 1.

\bigskip

\begin{center}
\scalebox{0.7}{\begin{tikzpicture}
\draw (-3.2,1.7) rectangle (3.2,-1.7);
\draw (-0.8,0) ellipse (1.5cm and 1.2cm);
\draw (0.8,0) ellipse (1.5cm and 1.2cm);
\end{tikzpicture}
\begin{picture}(0,0)(0,0)
\put(-38,68){\PencilLeftDown}
\put(-95,57){\rotatebox{0}{\scalebox{1.6}{\KreisGrossSchwarz}}}
\put(-104,43){\rotatebox{30}{\scalebox{1.6}{\QuadGrossSchwarz}}}
\put(-96,36){\rotatebox{-30}{\scalebox{1.6}{\DreieckGrossSchwarz}}}
\put(-135,60){\rotatebox{0}{\scalebox{1.5}{\KreisKleinSchwarz}}}
\put(-143,48){\rotatebox{-30}{\scalebox{1.5}{\QuadKleinSchwarz}}}
\put(-136,30){\rotatebox{30}{\scalebox{1.5}{\DreieckKleinSchwarz}}}
\put(-70,71){\rotatebox{-50}{\scalebox{1.6}{\DreieckGrossWeiss}}}
\put(-54.5,58.5){\rotatebox{-25}{\scalebox{1.6}{\QuadGrossWeiss}}}
\put(-66,42){\rotatebox{10}{\scalebox{1.7}{\QuadGrossGrau}}}
\put(-71,28){\rotatebox{-50}{\scalebox{1.6}{\KreisGrossWeiss}}}
\put(-44.5,45){\rotatebox{-25}{\scalebox{1.6}{\KreisGrossGrau}}}
\put(-56,29){\rotatebox{20}{\scalebox{1.6}{\DreieckGrossGrau}}}
\put(-177,84){\rotatebox{-15}{\scalebox{1.6}{\DreieckKleinWeiss}}}
\put(-177.5,66){\rotatebox{25}{\scalebox{1.7}{\QuadKleinGrau}}}
\put(-174.5,52){\rotatebox{25}{\scalebox{1.6}{\KreisKleinWeiss}}}
\put(-176.2,42.5){\rotatebox{-50}{\scalebox{1.3}{\DreieckKleinGrau}}}
\put(-176.8,21.5){\rotatebox{25}{\scalebox{1.3}{\KreisKleinGrau}}}
\put(-176.8,8.5){\rotatebox{-10}{\scalebox{1.5}{\QuadKleinWeiss}}}
\end{picture}}
\vspace{3pt}

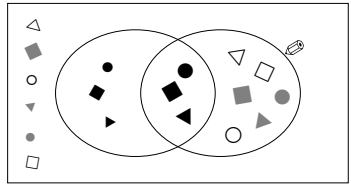
Figure 1: From shapes to Venn diagrams.
\end{center}

\bigskip

The examples in in Exercise 4(2)--4(4) can also be illustrated as in Figure 1 and should be treated as far as this is necessary. According to the participant's competencies, one can switch from Exercise 4 to Exercise 5. The latter can be completed using the shapes for illustration or by drawing diagrams with paper and pencil only.

\begin{exer} Draw \emph{Venn diagrams} for the sets in Exercise 4 und write down the sets, for instance\bigskip
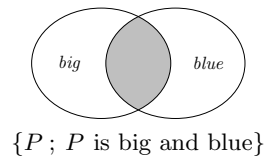
\begin{center}
\scalebox{0.61}{\begin{tikzpicture}
\begin{scope}
\clip (-0.8,0) ellipse (1.5cm and 1.2cm);
\clip (0.8,0) ellipse (1.5cm and 1.2cm);
\fill[color=lightgray] (-2,1.5) rectangle (3.3,-2);
\end{scope}
\draw (-0.5,0) (-0.8,0) ellipse (1.5cm and 1.2cm);
\draw  (0.8,0) ellipse (1.5cm and 1.2cm);
\draw(-1.5,0) node{\large big};
\draw(1.5,0) node{\large blue};
\end{tikzpicture}}
\\{\small $\{P\:;\:P \text{ \rm is big and blue}\}$}
\end{center}
in Exercise 4(1).
\end{exer}

Within a moderated discussion, the lecturer presents the solution for Exercise 5 and introduces the corresponding notations, cf.~Figure 2.

\medskip

\begin{center}
\scalebox{0.82}{
\fbox{\begin{minipage}[t]{521pt}
\begin{center}
\begin{tabular}{c|c|c}
\begin{minipage}{161pt}
\vspace{10pt}
\begin{center}
\begin{tikzpicture}
\begin{scope}
\clip (-0.5,0) ellipse (1.1cm and 0.9cm);
\clip (0.5,0) ellipse (1.1cm and 0.9cm);
\fill[color=lightgray] (-2,1.5)
rectangle (3.3,-2);
\end{scope}
\draw (-0.5,0) ellipse (1.1cm and 0.9cm);
\draw (0.5,0) ellipse (1.1cm and 0.9cm);
\draw(-1.1,0) node{\small big};
\draw(1.1,0) node{\small blue};
\end{tikzpicture}
\end{center}
\vspace{2pt}
\end{minipage}
&
\begin{minipage}{161pt}
\vspace{10pt}
\begin{center}
\begin{tikzpicture}
\begin{scope}
\clip (-0.5,0) ellipse (1.1cm and 0.9cm);
\fill[color=lightgray] (-2,1.5)
rectangle (3.3,-2);
\end{scope}
\begin{scope}
\clip (0.5,0) ellipse (1.1cm and 0.9cm);
\fill[color=lightgray] (-2,1.5)
rectangle (3.3,-2);
\end{scope}
\draw (-0.5,0) ellipse (1.1cm and 0.9cm);
\draw (0.5,0) ellipse (1.1cm and 0.9cm);
\draw(-1.1,0) node{\small plain};
\draw(1.1,0) node{\small red};
\end{tikzpicture}
\end{center}
\vspace{2pt}
\end{minipage}
&
\begin{minipage}{161pt}
\vspace{10pt}
\begin{center}
\begin{tikzpicture}
\begin{scope}
\clip (0,0) ellipse (1.5cm and 0.9cm);
\fill[color=lightgray] (-2,1.5)
rectangle (3.3,-2);
\clip (-0.3,-0.2) ellipse (0.7cm and 0.4cm);
\fill[color=white] (-2,1.5)
rectangle (3.3,-2);
\end{scope}
\draw (0,0) ellipse (1.5cm and 0.9cm);
\draw (-0.3,-0.2) ellipse (0.7cm and 0.4cm);
\draw(0.8,0.2) node{\small $U$};
\draw(-0.3,-0.2) node{\small abrasive};
\end{tikzpicture}
\end{center}
\vspace{2pt}
\end{minipage}
\\
\begin{minipage}[t]{161pt}\footnotesize
\begin{center}
$\{P\hspace{1pt};\hspace{0.5pt}P\text{ big and blue}\}$
\\$=\{P\hspace{1pt};\hspace{0.5pt}P\text{ big}\}\cap\{P\hspace{1pt};\hspace{0.5pt}P\text{ blue}\}$
\end{center}
\vspace{1pt}
\end{minipage}
&
\begin{minipage}[t]{161pt}\footnotesize
\begin{center}
$\{P\hspace{1pt};\hspace{0.5pt}P\text{ plain or red}\}$
\\$=\{P\hspace{1pt};\hspace{0.5pt}P\text{ glatt}\}\cup\{P\hspace{1pt};\hspace{0.5pt}P\text{ red}\}$
\end{center}
\vspace{1pt}
\end{minipage}
&
\begin{minipage}[t]{161pt}\footnotesize
\begin{center}
$\{P\hspace{1pt};\hspace{0.5pt}P \text{ not abrasive}\}$
\\$=U\hspace{0.5pt}\backslash\hspace{0.5pt}\{P\hspace{1pt};\hspace{0.5pt}P\text{ abrasive}\}=\{P\hspace{1pt};\hspace{0.5pt}P\text{ abrasive}\}^{\complement}$
\end{center}
\vspace{1pt}
\end{minipage}
\\\hline
\begin{minipage}[t]{161pt}\centering
\vspace{5pt}
\large Intersection\vspace{-1pt}
\\\phantom{a}
\end{minipage}
&
\begin{minipage}[t]{161pt}\centering
\vspace{5pt}
\large Union\vspace{-1pt}
\\\phantom{a}
\end{minipage}
&
\begin{minipage}[t]{161pt}\centering
\vspace{5pt}
\large Complement/Difference\vspace{-1pt}
\\\phantom{a}
\end{minipage}
\\\hline
\begin{minipage}{161pt}
\vspace{21pt}
\begin{center}
\begin{tikzpicture}
\draw (-1.1,0) ellipse (0.8cm and 0.6cm);
\draw (1.1,0) ellipse (0.8cm and 0.6cm);
\draw(-1.1,0) node{\small circular};
\draw(1.1,0.15) node{\small qua-};
\draw(1.1,-0.15) node{\small dratic};
\end{tikzpicture}
\end{center}
\vspace{2pt}
\end{minipage}
&
\begin{minipage}{161pt}
\vspace{15pt}
\begin{center}
\begin{tikzpicture}
\begin{scope}
\clip (-0.9,0.4) ellipse (0.6cm and 0.4cm);
\fill[color=lightgray] (-2,1.5)
rectangle (3.3,-2);
\end{scope}
\begin{scope}
\clip (0.9,0.4) ellipse (0.6cm and 0.4cm);
\fill[color=lightgray] (-2,1.5)
rectangle (3.3,-2);
\end{scope}
\begin{scope}
\clip (0,-0.4) ellipse (0.6cm and 0.4cm);
\fill[color=lightgray] (-2,1.5)
rectangle (3.3,-2);
\end{scope}
\draw (-0.9,0.4) ellipse (0.6cm and 0.4cm);
\draw (0.9,0.4) ellipse (0.6cm and 0.4cm);
\draw (0,-0.4) ellipse (0.6cm and 0.4cm);
\draw(-0.9,0.38) node{\small red};
\draw(0.95,0.38) node{\small\hspace{-2pt}yellow};
\draw(0,-0.4) node{\small blue};
\end{tikzpicture}
\end{center}
\vspace{2pt}
\end{minipage}
&
\begin{minipage}{161pt}
\vspace{12pt}
\begin{center}
\begin{tikzpicture}
\begin{scope}
\clip (-0.6,-0.2) ellipse (1.1cm and 0.9cm);
\fill[color=lightgray] (-2,1.5)
rectangle (3.3,-2);
\end{scope}
\begin{scope}
\clip (0.6,0.3) ellipse (1cm and 0.5cm);
\fill[color=white] (-2,1.5)
rectangle (3.3,-2);
\end{scope}
\draw (-0.6,-0.2) ellipse (1.1cm and 0.9cm);
\draw (0.6,0.3) ellipse (1cm and 0.5cm);
\draw(-1,0.10) node{\small circular};
\draw(-0.8,-0.30) node{\small and };
\draw(-0.55,-0.70) node{\small plain};
\draw(0.9,0.3) node{\small green};
\end{tikzpicture}
\end{center}
\vspace{2pt}
\end{minipage}
\\
\begin{minipage}[t]{161pt}\footnotesize
\begin{center}
$\{P\hspace{1pt};\hspace{0.5pt}P\text{ circular and quadratic}\}$
\\$=\{P\hspace{1pt};\hspace{0.5pt}P\text{ circular}\}\cap\{P\hspace{1pt};\hspace{0.5pt}P\text{ quadratic}\}=\emptyset$
\end{center}
\vspace{1pt}
\end{minipage}
&
\begin{minipage}[t]{161pt}\footnotesize
\begin{center}
$\{P\hspace{1pt};\hspace{0.5pt}P\text{ red or yellow or blue}\} $
\\$=\{P\hspace{1pt};\hspace{0.5pt}P\text{ red}\}\cup\{P\hspace{1pt};\hspace{0.5pt}P\text{ yellow}\}\cup\{P\hspace{1pt};\hspace{0.5pt}P\text{ blue}\}$
\end{center}
\vspace{1pt}
\end{minipage}
&
\begin{minipage}[t]{161pt}\footnotesize
\begin{center}
$\{P\hspace{1pt};\hspace{0.5pt}P\text{ circular and plain and not green}\}$
\\$=\{P\hspace{1pt};\hspace{0.5pt}P\text{ circular and plain}\}\hspace{0.5pt}\backslash\hspace{0.5pt}\{P\hspace{1pt};\hspace{0.5pt}P\text{ green}\}$
\end{center}
\vspace{1pt}
\end{minipage}
\\\hline
\begin{minipage}[t]{18pt}
\vspace{4pt}\hspace{-5pt}
{\Large\dbend}
\end{minipage}
\begin{minipage}[t]{130pt}
\vspace{4pt}
The intersection can be empty; the \textit{empty set} is denoted by $\emptyset$. Sets, whose intersection is empty, are said to be \textit{disjoint}.
\vspace{7pt} 
\end{minipage}
&
\begin{minipage}[t]{18pt}
\vspace{4pt}\hspace{-5pt}
{\Large\dbend}
\end{minipage}
\begin{minipage}[t]{130pt}
\vspace{4pt}
One can form the union of sets which are disjoint or not disjoint---no element appears twice.
\vspace{7pt} 
\end{minipage}
&
\begin{minipage}[t]{18pt}
\vspace{4pt}\hspace{-5pt}
{\Large\dbend}
\end{minipage}
\begin{minipage}[t]{130pt}
\vspace{4pt}
Complements are special differences. For $U\subseteq A$ we have $A^{\complement}=U\backslash A$.
\vspace{7pt} 
\end{minipage}
\\
\end{tabular}
\end{center}
\end{minipage}}}
\vspace{5pt}
\\Figure 2: Tabular corresponding to Exercise 5.
\end{center}

\subsection*{IV. Rules for Computation} The lecturer presents the following collection of rules at the blackboard. In a school framework it is unusual to select a presentation like those in Theorem 6. However, the latter is part of our concept, cf.~Section 1.

\begin{thm} For $A$, $B$, $C\subseteq U$ we have
\begin{compactitem}\vspace{2pt}
\item[(1)] $ A\cup A = A$, $A\cap A = A$,\vspace{2pt}
\item[(2)] $A\cup B = B\cup A$, $A\cap B = B\cap A$, \vspace{2pt}
\item[(3)] $A\cup (B\cup C) = (A\cup B)\cup C$, $A\cap (B\cap C) = (A\cap B)\cap C$,\vspace{2pt}
\item[(4)] $A\cup (B\cap C) = (A\cup B)\cap(A\cup C)$, $A\cap (B\cup C) = (A\cap B)\cup(A\cap C)$,\vspace{2pt}
\item[(5)] $A\cup\emptyset=A$, $A\cap\emptyset=\emptyset$,\vspace{2pt}
\item[(6)] $A\cup U = U$, $A\cap U=A$,\vspace{2pt}
\item[(7)] $(A^{\scriptscriptstyle\complement})^{\scriptscriptstyle\complement}=A$, $U^{\scriptscriptstyle\complement}=\emptyset$,\vspace{2pt}
\item[(8)] $(A\cup B)^{\scriptscriptstyle\complement}=A^{\scriptscriptstyle\complement}\cap B^{\scriptscriptstyle\complement}$, $(A\cap B)^{\scriptscriptstyle\complement}=A^{\scriptscriptstyle\complement}\cup B^{\scriptscriptstyle\complement}$.
\end{compactitem}
\end{thm}

Selected rules are examined in more detail individually or in small teams, for instance DeMorgan's rules in Theorem 6(8). While completing Exercise 7 below, the participants get used to the technique of colorizing diagrams, if necessary with further assistance by the lecturer. Shapes can be used for illustration but according to our experiences, cf.~Section 4, this is not necessary.

\begin{exer} Prove the equality $(A\cup B)^{\scriptscriptstyle\complement}=A^{\scriptscriptstyle\complement}\cap B^{\scriptscriptstyle\complement}$ via drawing the diagrams corresponding to the left and the right hand side consecutively.
\end{exer}

The solution consists of the stepwise construction of the desired diagrams, see Figure 3. Several of the steps shown below may be combined in one image by using multiple colors or patterns. We remark that the figure below should only give an idea how to solve the exercise by colorizing diagrams. Indeed, it is much easier to demonstrate this technique at the blackboard as it is to explain in printed form how it works. In addition we remark that it is maybe worth a discussion if colorizing diagrams as below can be considered as a proof. However, we propose to regard this again from a practical rather than from a philosophical point of view and therefore we suggest to skip this discussion in the context of the current workshop.

\medskip

\begin{center}

\scalebox{0.65}{\begin{minipage}[t]{650pt}
\begin{picture}(700,107)(0,0)
\put(0,0){\begin{tikzpicture}
\draw (-2.2,1.3) rectangle (-1.8,0.9);
\draw(-2.0,1.08) node{1.};
\draw (-2.2,1.3) rectangle (2.2,-1.3);
\begin{scope}
\clip (-0.5,0) ellipse (1.1cm and 0.9cm);
\fill[color=lightgray] (-2,1.5)
rectangle (3.3,-2);
\end{scope}
\draw (-0.5,0) ellipse (1.1cm and 0.9cm);
\draw (0.5,0) ellipse (1.1cm and 0.9cm);
\draw(-1.1,0) node{$A$};
\draw(1.1,0) node{$B$};
\draw(0,-1.8) node{\phantom{g}Colorize set $A$.\phantom{y}};
\draw(0,-2.2) node{\phantom{\,colorized in 1.~or in 2.}};
\end{tikzpicture}}
\put(130,0){\begin{tikzpicture}
\draw (-2.2,1.3) rectangle (-1.8,0.9);
\draw(-2.0,1.08) node{2.};
\draw (-2.2,1.3) rectangle (2.2,-1.3);
\begin{scope}
\clip (0.5,0) ellipse (1.1cm and 0.9cm);
\fill[color=lightgray] (-2,1.5)
rectangle (3.3,-2);
\end{scope}
\draw (-0.5,0) ellipse (1.1cm and 0.9cm);
\draw (0.5,0) ellipse (1.1cm and 0.9cm);
\draw(-1.1,0) node{$A$};
\draw(1.1,0) node{$B$};
\draw(0,-1.8) node{\phantom{g}Colorize set $B$.\phantom{y}};
\draw(0,-2.2) node{\phantom{\,colorized in 1.~or in 2.}};
\end{tikzpicture}}
\put(260,0){\begin{tikzpicture}
\draw (-2.2,1.3) rectangle (-1.8,0.9);
\draw(-2.0,1.08) node{3.};
\draw (-2.2,1.3) rectangle (2.2,-1.3);
\begin{scope}
\clip (-0.5,0) ellipse (1.1cm and 0.9cm);
\fill[color=lightgray] (-2,1.5)
rectangle (3.3,-2);
\end{scope}
\begin{scope}
\clip (0.5,0) ellipse (1.1cm and 0.9cm);
\fill[color=lightgray] (-2,1.5)
rectangle (3.3,-2);
\end{scope}
\draw (-0.5,0) ellipse (1.1cm and 0.9cm);
\draw (0.5,0) ellipse (1.1cm and 0.9cm);
\draw(-1.1,0) node{$A$};
\draw(1.1,0) node{$B$};
\draw(0,-1.8) node{Colorize everything that is};
\draw(0,-2.2) node{\,colorized in 1.~or in 2.};
\end{tikzpicture}}
\put(390,0){\begin{tikzpicture}
\begin{scope}
\clip (-2.2,1.3) rectangle (2.2,-1.3);
\fill[color=lightgray] (-2.2,1.5)
rectangle (2.2,-1.3);
\end{scope}
\begin{scope}
\clip (-0.5,0) ellipse (1.1cm and 0.9cm);
\fill[color=white] (-2,1.5)
rectangle (3.3,-2);
\end{scope}
\begin{scope}
\clip (0.5,0) ellipse (1.1cm and 0.9cm);
\fill[color=white] (-2,1.5)
rectangle (3.3,-2);
\end{scope}
\draw (-0.5,0) ellipse (1.1cm and 0.9cm);
\draw (0.5,0) ellipse (1.1cm and 0.9cm);
\draw(-1,0) node{$A$};
\draw(1,0) node{$B$};
\draw (-2.2,1.3) rectangle (2.2,-1.3);
\draw (-2.2,1.3) rectangle (-1.8,0.9);
\draw(-2.0,1.08) node{4.};
\draw(0,-1.8) node{Colorize everything that is};
\draw(0,-2.2) node{\,not colorized in 3.};
\end{tikzpicture}}
\end{picture}
\end{minipage}}

\vspace{5pt}

\scalebox{0.65}{\begin{minipage}[t]{650pt}
\begin{picture}(700,107)(0,0)
\put(0,0){\begin{tikzpicture}
\draw (-2.2,1.3) rectangle (2.2,-1.3);
\begin{scope}
\clip (-0.5,0) ellipse (1.1cm and 0.9cm);
\fill[color=lightgray] (-2,1.5)
rectangle (3.3,-2);
\end{scope}
\draw (-0.5,0) ellipse (1.1cm and 0.9cm);
\draw (0.5,0) ellipse (1.1cm and 0.9cm);
\draw(-1.1,0) node{$A$};
\draw(1.1,0) node{$B$};
\draw (-2.2,1.3) rectangle (-1.8,0.9);
\draw(-2.0,1.08) node{1.};
\draw(0,-1.8) node{\phantom{g}Colorize set $A$.\phantom{y}};
\draw(0,-2.2) node{\phantom{\,not colorized in 1.}};
\end{tikzpicture}}
\put(130,0){\begin{tikzpicture}
\begin{scope}
\clip (-2.2,1.3) rectangle (2.2,-1.3);
\fill[color=lightgray] (-2.2,1.5)
rectangle (2.2,-1.3);
\end{scope}
\begin{scope}
\clip (-0.5,0) ellipse (1.1cm and 0.9cm);
\fill[color=white] (-2,1.5)
rectangle (3.3,-2);
\end{scope}
\draw (-0.5,0) ellipse (1.1cm and 0.9cm);
\draw (0.5,0) ellipse (1.1cm and 0.9cm);
\draw(-1,0) node{$A$};
\draw(1,0) node{$B$};
\draw (-2.2,1.3) rectangle (2.2,-1.3);
\draw (-2.2,1.3) rectangle (-1.8,0.9);
\draw(-2.0,1.08) node{2.};
\draw(0,-1.8) node{Colorize everything that is};
\draw(0,-2.2) node{\,not colorized in 1.};
\end{tikzpicture}}
\put(260,0){\begin{tikzpicture}
\draw (-2.2,1.3) rectangle (2.2,-1.3);
\begin{scope}
\clip (0.5,0) ellipse (1.1cm and 0.9cm);
\fill[color=lightgray] (-2,1.5)
rectangle (3.3,-2);
\end{scope}
\draw (-0.5,0) ellipse (1.1cm and 0.9cm);
\draw (0.5,0) ellipse (1.1cm and 0.9cm);
\draw(-1.1,0) node{$A$};
\draw(1.1,0) node{$B$};
\draw (-2.2,1.3) rectangle (-1.8,0.9);
\draw(-2.0,1.08) node{3.};
\draw(0,-1.8) node{\phantom{g}Colorize set $B$.\phantom{y}};
\draw(0,-2.2) node{\phantom{\,not colorized in 1.}};
\end{tikzpicture}}
\put(390,0){\begin{tikzpicture}
\begin{scope}
\clip (-2.2,1.3) rectangle (2.2,-1.3);
\fill[color=lightgray] (-2.2,1.5)
rectangle (2.2,-1.3);
\end{scope}
\begin{scope}
\clip (0.5,0) ellipse (1.1cm and 0.9cm);
\fill[color=white] (-2,1.5)
rectangle (3.3,-2);
\end{scope}
\draw (-0.5,0) ellipse (1.1cm and 0.9cm);
\draw (0.5,0) ellipse (1.1cm and 0.9cm);
\draw(-1,0) node{$A$};
\draw(1,0) node{$B$};
\draw (-2.2,1.3) rectangle (2.2,-1.3);
\draw (-2.2,1.3) rectangle (-1.8,0.9);
\draw(-2.0,1.08) node{4.};
\draw(0,-1.8) node{Colorize everything that is};
\draw(0,-2.2) node{\,not colorized in 3.};
\end{tikzpicture}}
\put(520,0){\begin{tikzpicture}
\begin{scope}
\clip (-2.2,1.3) rectangle (2.2,-1.3);
\fill[color=lightgray] (-2.2,1.5)
rectangle (2.2,-1.3);
\end{scope}
\begin{scope}
\clip (-0.5,0) ellipse (1.1cm and 0.9cm);
\fill[color=white] (-2,1.5)
rectangle (3.3,-2);
\end{scope}
\begin{scope}
\clip (0.5,0) ellipse (1.1cm and 0.9cm);
\fill[color=white] (-2,1.5)
rectangle (3.3,-2);
\end{scope}
\draw (-0.5,0) ellipse (1.1cm and 0.9cm);
\draw (0.5,0) ellipse (1.1cm and 0.9cm);
\draw(-1,0) node{$A$};
\draw(1,0) node{$B$};
\draw (-2.2,1.3) rectangle (2.2,-1.3);
\draw (-2.2,1.3) rectangle (-1.8,0.9);
\draw(-2.0,1.08) node{5.};
\draw(0,-1.8) node{Colorize everything that is};
\draw(0,-2.2) node{\,colorized in 2.~and in 4.};
\end{tikzpicture}}
\end{picture}
\end{minipage}}

\vspace{3pt}

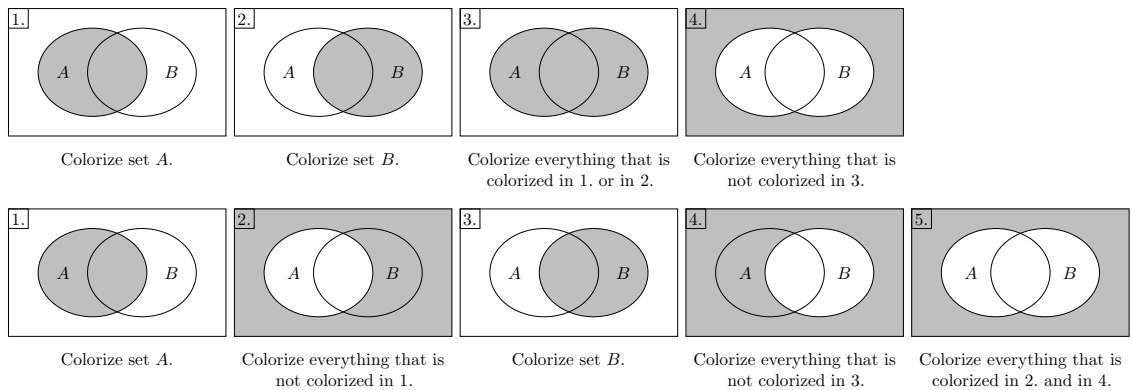
Figure 3: Idea for the solution of Exercise 7.
\end{center}

\subsection*{V. Further Discussion} The concept of sets as developed so far must now be consolidated by considering further aspects, for instance the paradoxes of naive set theory, e.g., the barber paradox or Russel's antinomy (see Smullyan \cite[p.~93]{Alice} and Potter \cite[Section 2.3]{Potter}) or by pointing out relations to other areas of mathematics. In any case, it must be emphasized that sets arise also in regular classes, see (a)--(d) at the end of this section.

\begin{ques} A certain barber living in a small town shaves all the inhabitants who don't shave themselves and never shaves any inhabitant who does shave himself. Does the barber shave himself or not?
\end{ques}

The particular formulation above is due to Smullyan \cite[p.~93]{Alice}. The question  can first be discussed in small groups, then via a moderated discussion. Finally, the lecturer puts the solution into writing.

\begin{sol} Denote by $M$ the set of all inhabitants and put
$$
S = \{m\in M\:;\: m \text{ shaves himself}\,\} \; \text{ and } \; B = \{ m\in M\:;\: m \text{ is shaved by the barber}\}.
$$
Then, $S\cup B=M$ and $B=M\backslash S$ hold. From the information in Question 8 we derive $B\cap S=\emptyset$. Let $b\in M$ be the barber. Assume that $b\in B$ holds. Then the barber shaves himself that is $b\in S$. A contradiction. Thus, $b\in S$ must hold. But then the barber is shaved by the barber. That is, $b\in B$, which is also a contradiction.
\end{sol}

In \textquotedblleft{}Solution\textquotedblright{}\,9 we defined the sets $S$ and $B$ by a \textquotedblleft{}self-reference\textquotedblright{}. $M$ is a set but $S$ and $B$ are not, since one cannot say if the barber belongs to $S$ or to $B$, respectively. Indeed, this not a matter of imprecise language. Russel's antinomy below is completely formal and displays the same kind of paradox. The lecturer explains the following and gives corresponding comments.

\begin{rem}\textbf{(Russel's antinomy, 1903)} Let $R$ be the set of all sets that are not members of themselves, i.e., $R=\{A\:;\:A\not\in A\}$. Assume, $R\in R$ holds. Then by definition $R\not\in R$ follows. Assume that $R\not\in R$ holds. Then by definition $R\in R$ follows. Thus we have $R\in R$ if and only if $R\not\in R$ and $R$ cannot be a set.
\smallskip
\\The naive definition of sets given at the beginning is adequate for many problems occurring in mathematics. However, Russel's antinomy highlights that the scope of the naive definition is limited. In \emph{axiomatic set theory} the notion of a set is defined by formal requirements which postulate the existence of certain sets, fix the properties of every set and determine the rules for constructing new sets of given ones. In this context it is then possible to prove formally that $R$ as defined above is not a set.
\end{rem}

In part V above we finally left the setting of the shapes. The definitions of union, intersection, complement etc., as well as the corresponding rules for computations remain valid also beyond the set $U$. The reason why we considered them only in this limited situation was just a matter of convenience. With Question 8 we even left the setting of naive set theory. However, we noted in Remark 10, that the latter is usable and useful in many situations. In order to illustrate this, the lecturer gives finally examples of sets which appear in regular lessons. Below we list some suggestions.

\begin{compactitem}\vspace{5pt}

\item[(a)] The set of natural numbers, the set of integers, the set of rational numbers, or the set of real numbers as well as subsets of the latter. For instance $\{n\in\mathbb{N}\:;\:n \text{ is even}\}=\{2n\:;\:n\in\mathbb{N}\}$, $[0,1]=\{x\in\mathbb{R}\:;\:0\leqslant{}x\leqslant{}1\}$ or $[0,1]\cap\mathbb{Q}=\{x\in\mathbb{Q}\:;\:0\leqslant{}x\leqslant1\}$. In all cases, the usage of mathematical notation is much more \emph{convenient} than the description of the sets in prose.

\vspace{3pt}

\item[(b)] A subset $f\subseteq\{(a,b)\:;\:a\in A,\,b\in B\}$ is a mapping from $A$ to $B$ provided that for every $a\in A$ there is one and only one $b\in B$ with $(a,b)\in f$. The value $f(a)$ of $f$ at $a$ is then the uniquely determined $b$ with $(a,b)\in f$. Usually one writes $f\colon A\rightarrow B$, $a\mapsto f(a)$ in this case. Sets allow for a formal definition of the concept of a map, which \emph{precises} the intuitive idea of the latter.

\vspace{3pt}

\item[(c)] Solution sets of systems of inequalities, e.g., $\{x\in\mathbb{R}\:;\:0\leqslant{}f(x)\leqslant{}g(x)\}=\{x\in\mathbb{R}\:;\:f(x)\geqslant0\}\cap\{x\:;\:f(x)\leqslant g(x)\}$. In order to determine the latter for given maps $f$, $g\colon\mathbb{R}\rightarrow\mathbb{R}$, the usage of sets allows for a reduction of the task in a \emph{systematic} way.

\vspace{3pt}

\item[(d)] The notion of similarity of triangles induces a system of subsets of the set of all triangles (two triangles belong to the same subset if and only if they are similar). The sets in this system are pairwise disjoint and their union is the set of all triangles. The language of sets, extended by the notions of \textquotedblleft{}equivalence relation\textquotedblright{}, \textquotedblleft{}equivalence class\textquotedblright{} and \textquotedblleft{}partition\textquotedblright{}, can make the geometric notion of similarity more \emph{accessible}.

\end{compactitem}\vspace{5pt}

It is left to the lecturer to decide which of the above examples he discusses in detail. This decision should depend on the background of the participants. Ideally, the examples refer to a current topic in regular classes or in another student workshop.
\smallskip
\\We remark that without pointing out the relation to other topics of mathematics the current workshop can also meet our objectives. However, the participants might get the impression that set theory is just a kind of mental exercise. We therefore strongly recommend to go into some of the points mentioned in (a)--(d) even if time limitations only allow for a very brief discussion.

\section{Practice Report}

In this section we report on two workshops which the author organized in 2006 and 2010 and whose topic was naive set theory. The outline given in Section 3 is based on these two workshops. We emphasize that all statements on execution and success that we make below are based on the personal assessment of the author; a formal evaluation was not carried out.
\smallskip
\\The first event on which we report was a one day workshop in which four students, two in the third year of primary education and two in the fourth year of primary education) participated. The workshop was part of an event organized by a regional group of the German Association for Highly Talented Children (DGhK-OWL). The workshop's content was similar to those explained in Section 3 but excluding Remark 10. The main differences are listed below.

\begin{compactitem}\vspace{5pt}

\item[1.] We did not fix a symbol for the set of all shapes, cf.~Definition 3. In Section 3 we added this, because it simplifies the definition for the complement (after Exercise 5).

\vspace{3pt}

\item[2.] In addition to the notions introduced in Definition 3, we introduced the logical operators $\wedge$, $\vee$ and $\neg$ in the sense of abbreviations for the words \textquotedblleft{}and\textquotedblright{}, \textquotedblleft{}or\textquotedblright{} and \textquotedblleft{}not\textquotedblright{}. In Exercise 5 we asked in addition to write down the sets by using the above operators. In connection with Theorem 6 and Exercise 7 it was then possible, to give also symbolic proofs of selected rules. In our outline above we dropped the points just mentioned in order to reduce the number of new symbols which have to be defined during the workshop.

\vspace{5pt}

\item[3.] Since the current event was a one day workshop, it was possible to treat a second topic in addition to naive set theory. For this second topic we selected the definition of maps and proceeded similar to (b) but with more details: We defined the cartesian product as the set of ordered pairs and then a map from $A$ to $B$ as a subset of $A\times B$ such that every element of $A$ occurs exactly in one pair on the left hand side. We practiced these notions by computing $\{1,2,3\}\times\{1,2\}$ and by discussing how many maps from $\{1,2\}$ to $\{a,b\}$ do exist. In addition, we introduced notions like \textquotedblleft{}value of a map at a point\textquotedblright{} or symbols like $f\colon A\rightarrow B$, $a\mapsto f(a)$, bijective maps, etc.

\end{compactitem}\vspace{5pt}

The second event on which this paper is based was a two hour workshop in the context of the Girl's Day 2010. The workshop took place at the University of Paderborn and was attended by 15 high school students in the seventh or eighth year of primary education. The author's aim was in particular to give the participants a realistic insight into mathematical studies---in addition to other Girl's Day events like talks on scientific topics or job perspectives for mathematicians in industry and business. Again, the workshop's content was similar to those explained in Section 3 and we thus only list the main differences.

\begin{compactitem}\vspace{5pt}

\item[4.] As in the first workshop discussed above, we introduced the operators $\wedge$, $\vee$ and $\neg$, see point 2. During the exercises some participants mixed up the logical operators with the set theoretical operators, a fact which supports our decision to drop the logical operators from the workshop concept as we did in Section 3.

\vspace{3pt}

\item[5.] In view of the age of the participants during this workshop, the shapes were only used until the beginning of Exercise 4. After the concept of Venn diagrams was introduced, only the latter and symbolic notation was used. In view of the latter we propose our concept of Section 3 for students in the seventh year of primary education or close to this.

\vspace{3pt}

\item[6.] The barber paradox was used as a discussion starter at the beginning of the workshop. The \textquotedblleft{}solution\textquotedblright{} was then given at the end. Due to time limitations the workshop did not include one of the points (a)--(d). 

\end{compactitem}\vspace{5pt}

In our opinion, both workshops met our objectives of Section 1. Those points which nevertheless needed improvement are mentioned above and are now part of our concept outlined in Section 3. The participants of both workshops liked the topic and contributed actively in the exercise sections. During the first workshop the students appreciated the interplay between concrete instances of sets (via shapes) and abstract thinking. During both workshops also the strongest participants at least punctually stretched to their limits. This is reasonable in view of the amount of new concepts, the style of presentation and the profoundity of the ideas discussed in particular at the end. However, from the author's point of view this is also necessary for a successful student workshop.

\bigskip

\footnotesize

{\sc Acknowledgements. } The author likes to thank H.~Berster, who initiated the first workshop mentioned in Section 4, as well as K.~Kelle and T.~Rupp who co-organized the second workshop. The author likes to thank K.~Hellermann, who supervised the thesis \cite{Wegner} on which this paper is based. In addition, the author likes to thank K.~Rolka and K.~Volkert for many useful comments. The final version of this paper was prepared during a stay at the Sobolov Institute of Mathematics at Novosibirsk; this stay was supported by the German Academic Exchange Service (DAAD).

\normalsize

\bibliographystyle{amsplain}

\bibliography{Workshop_Naive_Set_Theory}
\end{document}